\documentclass[11pt,a4paper]{article}
\usepackage{amsmath,amssymb,amsfonts,amsthm}
\usepackage[pdftex]{graphicx, color}

\newcommand{\E}{{\bf{E}}}
\newcommand{\PP}{{\bf{P}}}
\newcommand{\Var}{{\bf{Var}}}
\newcommand{\Bin}{{\bf{Bin}}}

\begin{document}

\bibliographystyle{plain}
\parindent=0pt
\centerline{\LARGE \bfseries A note on the vertex degree distribution} 
\centerline{\LARGE \bfseries of random intersection graphs}

\par\vskip 3.5em

\centerline{Mindaugas  Bloznelis}

\vglue2truecm

\centerline{Vilnius University, Institute of Computer Science,} 
\centerline{Naugarduko 24, Vilnius,
03225, Lithuania}


\bigskip





\begin{abstract}
We establish the
asymptotic degree distribution of the typical vertex of  
inhomogeneous and passive random intersection graphs 
under the minimal moment conditions.

\par
\end{abstract}

\smallskip
{\bfseries key words:}  degree distribution, random graph, random intersection graph, power law
\par\vskip 2.5em


\par\vskip 2.5em


{\bf 1. Introduction.}
Random intersection graphs  introduced by Karo\'nski et al \cite{karonski1999}
have attracted considerable attention in recent literature. They provide mathematically tractable theoretical models of complex networks that  capture important features of real networks: the power-law degree distribution, small typical distances between vertices, and a high statistical dependency of neighboring adjacency relations expressed in terms of 
non-vanishing clustering and assortativity coefficients, see \cite {BloznelisGJKR2015},
\cite{FriezeKaronski2016},  \cite{Spirakis2013} and references therein.

Vertex degree distributions have been studied by several authors
\cite{Bloznelis2008},
\cite{Bloznelis2013AAP},
\cite{BloznelisDamarackas2013},
\cite{DeijfenKets},
\cite{Stark2004},
\cite{JaworskiKaronskiStark2006},
\cite{JaworskiStark2008},
\cite{Rybarczyk2012},
\cite{Shang2010}.
In this note we establish the asymptotic degree distribution in two random intersection graph (RIG) models, the inhomogenious RIG and passive RIG, under the {\it optimal}  conditions.  Earlier papers that address
these random graph models  \cite{Bloznelis2013AAP},
\cite{BloznelisDamarackas2013} assumed excessive moment conditions. Here we present a simply and elegant proof that relax the moment conditions. 

\par
{\bf 1. Inhomogeneous random intersection graph}, \cite{Shang2010}. 
Let $X_1, X_2,\dots$ and $Y_1,Y_2,\dots $ be independent non-negative random variables such that
each $X_i$ has the  probability distribution $P_1$ and each $Y_j$ has the probability 
distribution
$P_2$. Given realized values $X=\{X_i\}_{i=1}^m$ and $Y=\{Y_j\}_{j=1}^n$
define the random bipartite graph $H_{X,Y}$ with the bipartition 
$V=\{v_1,\dots, v_n\}$, $W=\{w_1,\dots, w_m\}$, 
where  edges $\{w_i,v_j\}$ are inserted 
with probabilities
$p_{ij}=\min\{1, X_iY_j(nm)^{-1/2}\}$ independently 
for each $\{i,j\}\in [m]\times [n]$.
The 
 inhomogeneous random intersection graph  $G(P_1,P_2, n,m)$ defines the adjacency 
relation on the  vertex 
set $V$: vertices $v',v''\in V$ are declared adjacent  whenever
$v'$ and $v''$ have a common neighbour in $H_{X,Y}$. 
Let  $d=d(v_1)$ denote the degree of vertex $v_1$  in
$G(P_1,P_2, n,m)$.  
Denote $a_k=\E X_1^k$, $b_k=\E Y_1^k$.
 The following result is shown in Theorem 1 (ii) of
 {\cite{BloznelisDamarackas2013}}. 

{\bf Theorem A}.
{\it 
Let $\beta\in (0,+\infty)$.
 Let $m,n\to\infty$. Assume that $m/n\to \beta$. 
 Suppose  that $\E X_1^2<\infty$ and $\E Y_1<\infty$. 
Then  $d$ converges in distribution to the random variable
\begin{equation}\label{d*1}
d_*=\sum_{j=1}^{\Lambda_1}\tau_j, 
\end{equation}
where $\tau_1,\tau_2,\dots$ are independent 
and identically distributed random variables independent of the random variable $\Lambda_1$.
They are distributed as follows. For $r=0,1,2,\dots$, we have
\begin{equation}\label{d*1++}
\PP(\tau_1=r)=\frac{r+1}{\E\Lambda_2}\PP(\Lambda_2=r+1)
\qquad
{\text{and}}
\qquad
\PP(\Lambda_i=r)=\E \,e^{-\lambda_i}\frac{\lambda_i^r}{r!},
\qquad
i=1,2.
\end{equation}
Here $\lambda_1=Y_1a_1\beta^{1/2}$ and $\lambda_2=X_1b_1\beta^{-1/2}$.
}

\smallskip

It was conjectured in {\cite{BloznelisDamarackas2013}} 
that the second moment condition $\E X_1^2<\infty$ 
can be relaxed to the first moment condition 
$\E X_1<\infty$. We show that this is the case.

{\bf Theorem 1}.
{\it
Theorem A remains true if we replace the second moment condition
 $\E X_1^2<\infty$ by the first moment condition $\E X_1<\infty$.
 }

\bigskip

{\bf 2. Passive random intersection graph}, \cite{GodehardtJaworski2001}.
 Let $P$ be a probability distribution on $\{0,1,\dots, m\}$.
Let $D_1,\dots, D_n$ be independent random subsets of $W=\{w_1,\dots, w_m\}$ having the same
probability distribution $\PP(D_i=A)={\binom{m}{|A|}}^{-1}P(|A|)$, $A\subset W$.
The  passive random 
intersection graph $G^*(n,m,P)$ defines the adjacency relation on the vertex set $W$:
two vertices $w, w'\in W$ are declared adjacent whenever 
 $w,w'\in D_j$ for some
$j$ (\cite{GodehardtJaworski2001}, \cite{JaworskiStark2008}). Let  $d=d(w_1)$ denote the degree of 
vertex $w_1$ in  $G^*(n,m,P)$. 
By $P_\xi$ we denote the probability 
distribution of a random variable $\xi$. Thus, 
$P_{X_1}=P$ for
$X_1:=|D_1|$. Furthermore, 
given a probability distribution $Q$ on $\{0,1,2,\dots\}$ with  a 
finite first moment $\mu_Q=\sum_iiQ(i)<\infty$, let
${\tilde Q}$ denote the size biased distribution, ${\tilde Q}(j)=(j+1)Q(j+1)\mu_Q^{-1}$, $j=0,1,\dots$.
The following result is shown in Theorem 3.1 of  {\cite{Bloznelis2013AAP}}.

{\bf Theorem B.} 
{\it
Let $\beta\in (0,+\infty)$.  Let $m,n\to\infty$. Assume that $m/n\to \beta$ and

(i) $X_{1}$ converges in distribution to a random variable $Z$;

(ii)  $0<\E Z<\infty$ and \ $\lim_{m\to\infty}\E X_{1}=\E Z$;

(iii) $\E Z^{4/3}<\infty$ and \ $\lim_{m\to\infty}\E X^{4/3}_{1}=\E Z^{4/3}$.

Then $d$ converges in distribution to the compound Poisson random variable 
$d_*=\sum_{j=1}^\Lambda {\tilde Z}_j$.
Here ${\tilde Z}_1,{\tilde Z}_2,\dots$ are independent random variables with 
the common probability distribution ${\tilde P}_Z$, the random variable $\Lambda$ 
is independent of the sequence 
${\tilde Z}_1,{\tilde Z}_2,\dots$ and has Poisson distribution 
with mean $\E \Lambda= \beta^{-1}\E Z$. 
}

It was conjectured in {\cite{Bloznelis2013AAP}} that  condition (iii) is redundant.
 We show that this is the case.

{\bf Theorem 2}.
{\it
The conclusion of 
Theorem B  remains true if we drop condition (iii).
 }

\bigskip
 
{\bf 3. Proof of Theorems 1 and 2}. 
Let $M>0$ be an integer. Given $n,m$, let
 ${\hat G}$ and 
${\check G}$ denote the inhomogeneous (passive) intersection graph with $X_i$,
$1\le i\le m$, replaced by  
${\hat X}_i=X_i{\mathbb I}_{\{X_i\le M\}}$ 
and 
${\check X}_i=X_i{\mathbb I}_{\{X_i> M\}}$
 respectively.
Let ${\hat d}$ and ${\check d}$ denote the degree 
of $v_1$ ($w_1$)  in
 ${\hat G}$ and 
${\check G}$  respectively.
We have ${\hat d}\le d\le {\hat d}+{\check d}$.
Hence, for any integer $k\ge 0$ we have
\begin{equation}\label{2019-06-17}
 \PP({\hat d}\ge k)
 \le
 \PP(d\ge k)
 \le 
 \PP({\hat d}\ge k)
 +
 \PP({\check d}\ge 1). 
\end{equation}
In the case of inhomogeneous graph we have
\begin{eqnarray}\label{2019-06-17+1}
&&
\PP({\check d}\ge 1)
=
\E{\mathbb I}_{\{{\check d}\ge 1\}}
\le 
\E
\sum_{i=1}^m
{\mathbb I}_{\{\{v_1, w_i\}\in H_{X,Y}\}}
{\mathbb I}_{\{|X_i|>M\}}
\le 
m
\E \frac{Y_1{\check X}_1}{\sqrt{nm}}
=
b_1
\frac
{\sqrt{m}}
{\sqrt {n}}
\E {\check X}_1.
\end{eqnarray}
In the case of passive graph we have
 \begin{eqnarray}\label{2019-06-17+2}
&&
\PP({\check d}\ge 1)
=
\E{\mathbb I}_{\{{\check d}\ge 1\}}
\le 
\E
\sum_{i=1}^n
{\mathbb I}_{\{w_1\in D_i\}}
{\mathbb I}_{\{|D_i|>M\}}
=
\frac{n}{m}
\E 
{\check X}_1.
\end{eqnarray}
Note that conditions of Theorem 1  (Theorem 2) imply 
$\E {\check X}_1=o(1)$ for $M\to+\infty$.
 
Let ${\hat d}_*={\hat d}_*(M)$ be the limiting distribution of ${\hat d}$ defined by Theorem A (Theorem B). We obtain  from  Theorem A (Theorem B) and 
(\ref{2019-06-17}) that for any integer 
$k\ge 0$ 
\begin{equation}\label{2019-06-14+2}
\PP({\hat d}_*\ge k)
\le 
\liminf_{n,m\to+\infty}
\PP(d\ge k)
\le 
\limsup_{n,m\to+\infty}
\PP(d\ge k)
\le 
\PP({\hat d}_*\ge k)
+
\limsup_{n,m\to+\infty} \PP({\check d}\ge 1).
\end{equation}
Letting $M\to+\infty$ we obtain
$\PP({\hat d}_*\ge k)\to \PP(d_*\ge k)$ and  
$\limsup_{n,m\to+\infty} \PP({\check d}\ge 1)=0$, see
(\ref{2019-06-17+1}),
(\ref{2019-06-17+2}).
Now (\ref{2019-06-14+2}) implies that
$d$ converges in distribution to $d_*$.

\end{document}